\documentclass[12pt]{article}
\usepackage{amsmath,amsxtra,amsfonts,amssymb,amscd,latexsym,pstcol,pst-plot,pst-fill,pst-grad,graphicx}
\usepackage{makeidx}

\newtheorem{proposition}{Proposition}

\newtheorem{theorem}{Theorem}

\newcommand{\reels}{\mathbb{R}}
\newcommand{\nat}{\mathbb{N}}
\newcommand{\rel}{\mathbb{Z}}

\newcommand{\esper}{\mathbb{E}}
\newcommand{\proba}{\mathbb{P}}

\begin{document}

\title{Strong Uniform Consistency of the Frequency Polygon Density estimator for Stable Non-Anticipative Stochastic Processes}
\author{Salim Lardjane}
\date{{\em salim.lardjane@univ-ubs.fr}\\
Laboratoire de Math\'{e}matiques de Bretagne Atlantique\\
 UMR CNRS 6205\\
Rue André Lwoff - BP 573\\ 
56017 Vannes, France}

\maketitle

\noindent {\bf Abstract}. {The author establishes a new mathematical expression for the Frequency Polygon. He uses it to prove the strong uniform consistency of the Frequency Polygon marginal density estimator for non-anticipative stationary stochastic processes which are stable in the sense of Wu \cite{Wu05}. He gives examples of several times series models for which this result is relevant.}\\


\section{Introduction}

The histogram is without contest the most used density estimation me\-thod to represent the distribution of a sample from a continuous random variable. This is due mainly to its conceptual simplicity and its ease of computation. However, it has the great disadvantage of being discontinuous and less than optimal from the theoretical point of view \cite{Scott85}. To solve these problems, other density estimation methods are used, among which the most well-known is undoubtedly Parzen-Rosenblatt's Kernel Density Estimator (KDE) \cite{Parzen62,Rosenblatt58}. This estimator, although optimal from the theoretical point of view, has the disadvantage of being computationaly demanding. Namely, to estimate the density in $m$ points using this estimator requires $O(nm)$ computations, where $n$ is the size of the sample. This problem has become particularly acute the last years because of the appearance of increasingly massive databases one wants to explore in a quick way.\\
Thus, there is currently a movement towards a less known estimator, which is based on the histogram, defined in a simple way, continuous, computationaly less demanding than the KDE and which shares the good theoretical properties of the latter : it is the Frequency Polygon \cite{Scott85}.\\
The histogram is usually defined by dividing the real axis in a succession of bins of equal width and counting the number of observations in each bin. More precisely, the value of the histogram at a point is given by the relative frequency of the interval containing the point. The Frequency Polygon is defined by joining with a segment the values of the histogram at the middle of adjacent bins.\\
This estimator is continuous by construction and requires $O(mp_n)$ operations to estimate the density in $m$ points, where $p_n$ is the number of bins with nonzero count for a sample of size $n$ (using a hash table for example).\\
Despite the fact that the Frequency Polygon is a rather old method of density representation, its asymptotic properties as a density estimator have been first investigated in the seventies of the twentieth century, in the case of a sample of $n$ independent identically distributed random variables. Thus, it has been established that the Frequency Polygon reached the same optimal rate of uniform convergence as the KDE \cite{Revesz72}.\\
The investigation of the uniform consistency of the Frequency Polygon for dependent data started in the nineties, with a paper by Carbon, Garel and Tran (\cite{Carbon97} in the case of a sample from an $\alpha$-mixing stochastic process. More recently, Xing \cite{Xing15}, Yang \cite{Yang15}, Kang et al. \cite{Kang18} and Wang et al. \cite{Wang21} have extended these results, either by deriving them under weaker conditions or by using other mixing assumptions.\\
In this paper, we adopt a third way. Namely, we give conditions that are very different from those existing in the litterature, under which the Frequency Polygon is uniformly consistent for the marginal density of stationary non-anticipative stochastic processes that are stable in the sens of Wu \cite{Wu05,Wu06}. This family of stochastic processes and our assumptions are sufficiently general to include many time series models, as we shall see in section \ref{exemples}.\\

\noindent {\bf Notations.} Let $(\zeta_n)$ be a sequence of random variables and $(a_n)$ a deterministic sequence. We shall note $\zeta_n=O_{p.s}(a_n)$ if $\zeta_n / a_n$ is almost surely bounded when $n\to\infty$ and $\zeta_n=o_{p.s}(a_n)$ if $\zeta_n / a_n\to 0$ almost surely when $n\to\infty$.

\section{Model and assumptions}

Let $\{\varepsilon_0^\prime, \varepsilon_n,n\in\rel\}$ be a sequence of independent identically distributed square integrable random variables with null expectancy and variance $\sigma^2>0$ on a probability space $(\Omega,\mathcal{A},\proba)$. Let us denote by
$\xi_k=(\dots,\varepsilon_{k-1},\varepsilon_k)$ the history of this process  up to time $k\in \rel$. We consider stochastic processes $\{X_n,n\in\rel \}$ that satisfy the relation :
\begin{equation}\label{eqdef}
X_n=\mu(\xi_{n-1})+\varepsilon_n
\end{equation}
where $\mu$ is a real-valued  measurable mapping on the Borel space of real sequences. Then, $\{(X_n,\xi_n),n\geq 1\}$ and $\{X_n,n\geq 1\}$ are strictly stationary and ergodic. Our first assumption writes :
\begin{itemize}
\item[$\mathcal{A}_1$:]
$ \esper[ \mu^2(\xi_{-1})]<\infty$.
\end{itemize}
This assumption entails  in particular that $Var(X_n)=  Var(\mu(\xi_{-1}))+\sigma^2 < \infty$ for every $n\in\rel$.\\
Now, for every $k\geq 0$, let us define :
$\xi_k^*=(\xi_{-1},\varepsilon_0^\prime,\varepsilon_1,\dots,\varepsilon_k)$ 
and set
$
X_0^*=\mu(\xi_{-1})+\varepsilon_0^\prime
$
and
$$
X_k^*=\mu(\xi_{k-1}^*)+\varepsilon_k \quad (k>1).
$$
Our dependence assumption will be expressed using the standard deviations
$
\delta_k=Var(X_k-X_k^*)^{1/2}
$. Namely, we shall assume :
\begin{itemize}
\item[$\mathcal{A}_2$:] $
\sum_{k=1}^\infty \delta_k <\infty.
$
\end{itemize}
Assumption $\mathcal{A}_2$ is a short-range dependence condition. Indeed, $\delta_k$ quantifies the variation in $X_k-X_k^*$ and thus the impact of $\varepsilon_0$ on the value of the series at time $k$. Thus, $\mathcal{A}_2$ can be interpreted by saying that the contribution of $\varepsilon_0$ to the future values of the time series is asymptotically negligible.\\
The quantities $\delta_k$ are denoted $\delta_2(k)$ and termed {\em physical or functional dependence measures} by Wu \cite{Wu05}, who proposes to call the processes satisfying $\mathcal{A}_2$ 2-stable.

Now, let us denote by $\mu_\varepsilon$ the distribution of $\varepsilon_1$, by $H^2(\reels)$ the Sobolev space of order 2 on $\reels$ and by $\rho_\beta(f)$ the Hölder seminorm :
$$
\rho_\beta(f)=\sup_{x,y\in\reels, x\neq y} |x-y|^{-\beta} |f(x)-f(y)|.
$$
where $0<\beta\leq 1$. Our third and last assumption is the following :
\begin{itemize}
\item[$\mathcal{A}_3$:]  $\mu_\varepsilon$ is absolutely continuous with respect to the Lebesgue measure $\lambda$ on $\reels$ and ${d\mu_\varepsilon}/{d\lambda}\in H^2(\reels)$.\\
\end{itemize}
Condition $\mathcal{A}_3$ implies that:

\noindent $\text{I.}$ $d\mu_\varepsilon/d\lambda$ has a representative $f_\varepsilon$ in $\mathcal{L}^2(\reels)$ such that, for every scalars $x,y$ :
$
f_\varepsilon(x)-f_\varepsilon(y)=\int_x^y f_\varepsilon^\prime(s) ds
$
where $f_\varepsilon^\prime$ denotes the weak derivative (the derivative in the sense of distributions) of $f_\varepsilon$.\\
$\text{II.}$  $f_\varepsilon$ is continuous and goes to 0 at infinity.\\
$\text{III.}$  $f_\varepsilon(x)$ is bounded, namely there exists a positive constant $c_0$ such that $\sup_x f_\varepsilon(x)\leq c_0$. \\
Then, the conditional density of $X_1$ given $\xi_0$, denoted by $ f_1(x|\xi_0)$, is well defined and $\sup_x f_1(x|\xi_0)\leq c_0$ $\proba$-almost surely. More precisely, one has :
$$
 f_1(x|\xi_0=\xi)=f_\varepsilon(x-\mu(\xi)).
$$
Hence, the marginal density of $\{X_n,n\geq 1\}$, denoted by $f$, is well defined and satisfies :
$$
f(x)=\esper[f_1(x|\xi_0)]\leq c_0\quad\proba-p.s.
$$
$\text{IV.}$ $f_\varepsilon$ is Lipschitz-regular, that is
$
\rho_{1}(f_\varepsilon)<\infty.
$
This entails that $f$ is Lipschitz-regular and hence uniformly continuous on $\reels$.  To see this, write :
\begin{eqnarray*}
|f(y)-f(x)| 
 & \leq & \esper [|f_\varepsilon(y-\mu(\xi_0))-f_\varepsilon(x-\mu(\xi_0))|]
\end{eqnarray*}
hence , for every scalars $x,y$:
\begin{equation}\label{holder}
|f(y)-f(x)| \leq \rho_{1}(f_\varepsilon)\cdot |y-x|
\end{equation}
Now, let us denote by $\varphi_\varepsilon(t)$ the characteristic function of $\varepsilon_1$, that is :
$
\varphi_\varepsilon(t)=\int_\reels e^{itx}\,f_\varepsilon(x)dx
$
where $i^2=-1$. Then, 
$$
\text{V.} \;\int_\reels |\varphi_\varepsilon(t)|^2 \,(1+t^2)\,t^2\,dt < \infty.
$$
This property is satisfied for instance if $|\varphi_\varepsilon(t)|=O(|t|^{-\alpha})$ when $|t|\to\infty$, where
$\alpha>5/2$.\\ 
Proofs of properties I-V are given in Taylor \cite{Taylor96} and Adams and Fournier \cite{Adams03} for instance. Assumption $\mathcal{A}_3$ is clearly satisfied if the noise is Gaussian and is satisfied in particular if $\mu_\varepsilon$ admits a density $f_\varepsilon$ which is twice continuously differentiable in the usual sense on $\reels$ and such that $df_\varepsilon(x)/dx$ and $d^2f_\varepsilon(x)/dx^2$ are in $\mathcal{L}^2(\reels)$.

\section{The Frequency Polygon}

Suppose that we observe the time series $\{X_n,n\in\rel\}$ from the time zero to time $n$ and let $(k_n)_{n\in\rel}$ be a sequence of regularly distributed points on the real axis, such that $k_{n+1}-k_n=b_n>0$ for every $n\in\rel$, with $k_0= 0$.  Then, $k_n=nb_n$. 

The histogram density estimator of $f$ is defined, for every scalar $x$, by the expression :
$$
f_n(x)=\frac{\nu_n(x)}{nb_n}
$$
where $\nu_n(x)$ is the number of observations in the same bin as $x$, or equivalently by:
$$
f_n(x)=\frac{1}{nb_n}\sum_{i=1}^n \mathbf{1}_{]k_n(x),k_n(x)+b_n]}(X_i)
$$
where $\mathbf{1}_E$ denotes the indicator function of the set $E$ and where $k_n(x)$ is the greatest bin limit that is strictly lower than $x$, that is :
$
k_n(x)=\lfloor x/b_n\rfloor\cdot b_n.
$

The Frequency Polygon density estimator, denoted by $g_n(x)$ is obtained by joining by a segment the values of the histogram at the middle of the bins.\\
The usual formalization of this definition is the following :
$$
g_n(x)= \left(\frac{1}{2}+k-\frac{x}{b_n}\right)\,f_n(kb_n)+\left(\frac{1}{2}-k+\frac{x}{b_n}\right)\,f_n((k+1)b_n)
$$
for every $x\in ]kb_n-\frac{b_n}{2},kb_n+\frac{b_n}{2}]$ and  $k\in\rel$.
In the sequel, we shall introduce a new expression of the Frequency Polygon, which appears to be of independent interest. To this end, let us start by noting that
$$
\frac{\nu_n(x)}{n}=F_n(k_n(x)+b_n)-F_n(k_n(x))
$$
where $F_n(x)$ is the empirical distribution function, defined for every scalar $x$, by the expression:
$
F_n(x)=\frac{1}{n}\sum_{i=1}^n \mathbf{1}_{]-\infty,x]}(X_i).
$

For every $n>0$, we can define an operator $A_n$ on the linear space of bounded real functions of a real argument by
$$
A_{n}F=\sum_{k\in\rel}\frac{F(kb_n+b_n)-F(kb_n)}{b_n}\cdot \mathbf{1}_{]kb_n,kb_n+b_n]}
$$
Then, we have $f_n=A_ nF_n$, and we can state our first result, which provides a new expression for the Frequency Polygon. Let us denote by $\tau_a$ the shift operator with parameter $a$, that is $\tau_af(x)=f(x-a)$. 
\begin{theorem}\label{FP}
The Frequency Polygon $g_n(x)$ can be defined using the operator
$$
B_n =(1-u_n)\cdot \tau_{b_n/2}A_n  +u_n\cdot \tau_{-b_n/2} A_n 
$$
where :
$$
u_n(x)=x/b_n-\frac{\lfloor x/b_n-1/2\rfloor+\lfloor x/b_n+1/2\rfloor}{2}
$$
for every $x\in\reels$. Namely,
$
g_n=B_nF_n.
$
\end{theorem}
{\em \bf Proof. } Suppose that $x$ is such that $kb_n-b_n/2<x\leq k b_n+b_n/2$. Then,
$
\lfloor x/b_n-1/2\rfloor = k-1
$
and
$
\lfloor x/b_n+1/2\rfloor =k.
$
Consequently,
$
u_n(x)=1/2-k+x/b_n
$
and
$
1-u_n(x)=1/2+k-x/b_n.
$
Moreover,
$$
\tau_{b_n/2}A_nF_n(x)=\tau_{b_n/2}f_n(x)=f_n(x-b_n/2)=f_n(kb_n)
$$
and
$$
\tau_{-b_n/2}A_nF_n(x)=\tau_{-b_n/2}f_n(x)=f_n(x+b_n/2)=f_n((k+1)b_n)
$$
hence,
$$
g_n(x)=(1/2+k-x/b_n)\cdot f_n(kb_n)+(1/2-k+x/b_n)\cdot f_n((k+1)b_n).
$$
Thus, we obtain the standard expression of the Frequency Polygon, which establishes the result. $\square$\\
The following proposition lists some properties of the operators $A_n$ and $B_n$.
\begin{proposition}\label{lemma}
Let $(b_n)_{n\in\nat}$ be a positive sequence converging to 0, assume $\mathcal{A}_3$ is satisfied, and let $F$ denote the marginal cumulative distribution function of $\{X_n,n\in\rel\}$. For every $n>0$, one has:
\begin{itemize}
\item[$(i)$] $A_n$ is a bounded linear operator such that $||A_n||_\infty=O( b_n^{-1})$, and $||A_{n} F-f||_\infty= O(b_n)$ when $n\to\infty$. 
\item[$(ii)$]  $B_n$ is a bounded linear operator such that $||B_n||_\infty=O(b_n^{-1})$, and $||B_{n} F-f||_\infty= O(b_n)$ when $n\to\infty$. 
\end{itemize}
\end{proposition}
\noindent {\bf Proof}. $A_n$ is clearly linear. Moreover, for every $n>0$,  one has :
$$
||A_nF||_\infty\leq \frac{1}{b_n}\sup_{z\in\rel} |F(zb_n+b_n)-F(zb_n)|\leq b_n^{-1}||F||_\infty.
$$
Thus, the operator $A_n$ is bounded, with $||A_n||_\infty \leq  b_n^{-1}$.\\
On the other side, we can write :
\begin{eqnarray*}
||A_nF-f||_\infty & = &  \sup_{x\in\reels} \left| \frac{F(k_n(x)+b_n)-F(k_n(x))}{b_n}-f(x)\right|\\
 & = & \sup_{x\in\reels} | f(\xi_n(x))-f(x)|\leq \rho_{1}(f)\cdot b_n
\end{eqnarray*}
where $\xi_n(x)\in]k_n(x),k_n(x)+b_n[$. This establishes the second property of $A_n$ $(n>0)$. \\
Now, $B_n$ is clearly linear and we can write :
\begin{eqnarray*}
||B_n||_\infty &\leq & | (||\tau_{b_n/2}||_\infty +||\tau_{-b_n/2}||_\infty) ||A_{n}||_\infty\\
 &\leq & (||\tau_{b_n/2}||_\infty +||\tau_{-b_n/2}||_\infty) \,b_n^{-1}  =  O( b_n^{-1})
\end{eqnarray*}
which shows the first part of $(ii)$. Moreover, we have :
\begin{eqnarray*}
||B_nF-f||_\infty & \leq & || \tau_{b_n/2}A_nF\cdot (1-u_n)+\tau_{-b_n/2} A_n F\cdot u_b -f||_\infty\\
 & \leq  &  ||  (1-u_n)\cdot (\tau_{b_n/2} A_nF-f)+u_n\cdot (\tau_{-b_n/2} A_nF-f)||_\infty\\
  &=& ||  (1-u_n)\cdot \tau_{b_n/2} (A_nF-f)+(1-u_n) (\tau_{b_n/2}f-f)\\
  &  & \quad +u_n\cdot \tau_{-b_n/2} (A_nF-f)+u_b\cdot (\tau_{-b_n/2} f-f)||_\infty\\
  & \leq & 2||A_nF-f||_\infty+||\tau_{b_n/2}f-f||_\infty+||\tau_{-b_n/2}f-f||_\infty\\
  & = & O(b_n)
\end{eqnarray*}
hence the second part of property $(ii)$. $\square$

\section{Strong uniform consistency}

To investigate the uniform consistency of the Frequency Polygon $g_n$, we shall bound the error $||g_n-f||_\infty=||B_nF_n-f||_\infty$ from above by the sum of two terms : $||B_nF_n-B_nF||_\infty$ and $||B_nF-f||_\infty$. \\
The first term is small if $F_n$ is close to $F$, which happens if one has enough data in each bin, that is if the bins are wide enough, but the second term is small if $b_n$ is small, that is if the bins are narrow enough. The problem is thus to find the appropriate relation between the width of the bins $b_n$ and the sample size $n$.\\
To investigate the question, we shall need the following theorem, which is due to Wu \cite{Wu06}.
\begin{theorem}\label{lemme1}
Assume $\mathcal{A}_1-\mathcal{A}_3$. Let $F$ denote the marginal cumulative distribution function of $\{X_n,n\in\rel\}$ and $F_n$ the empirical cumulative distribution function of the sample $X_1,\dots,X_n$. Let $G_n=\sqrt{n}(F_n-F)$ and denote by
$$
 \Delta_n (b) =\sup_{|u-v|\leq b} |G_n(u)-G_n(v)|
$$
the modulus of continuity of $G$, for $b>0$.\\
Assume that $b_n\to 0$, $\log n =O(nb_n)$ and set:
$
\kappa_n=(\log n)^{1/2}\log\log n.
$
Then,
$$
\Delta_n(b_n)=O_{a.s.}(\sqrt{b_n\log n})+o_{a.s.}(b_n\kappa_n)
$$
\end{theorem}
\noindent {\bf Proof}. See Th. 2.1. and 2.2. in \cite{Wu06}. $\square$\\
We are now in position to state our main result.
\begin{theorem}\label{convunif}
If assumptions $\mathcal{A}_1-\mathcal{A}_3$ are satisfied and if $b_n\to 0$ and $\frac{n b_n}{\log n}\to \infty$ when $n\to\infty$,  then $g_n$ converges uniformly almost surely to $f$. Moreover, if $b_n=(n^{-1}\log n)^{1/3}$, then
$$
||g_n-f||_\infty = O_{p.s.}\left((n^{-1}\log n)^{1/3}\right).
$$
\end{theorem}
Let us note that the rate $(n^{-1}\log n)^{1/3}$
is the optimal rate obtained by Stone \cite{Stone83} in the setting of i.i.d. observations under our regularity assumptions.\\
\noindent {\bf Proof}.  Let $(b_n)_{n\in\nat}$ be a positive sequence converging to 0. For every $n$, we have:
$$
||A_{n}F_n-f||_\infty \leq ||A_{n}F_n-A_{n}F||_\infty +||A_{n}F-f||_\infty
$$
Set $G_n=\sqrt{n}(F_n-F)$. Then, we can write 
 $$
 A_{n}F_n-A_{n}F=\frac{G_n(nb_n+b_n)-G_n(nb_n)}{\sqrt{n}b_n}
 $$ 
hence
 $$
 ||A_{n}F_n-A_{n}F||_\infty \leq \frac{1}{\sqrt{n}b_n}\sup_{|u-v|\leq b_n} |G_n(u)-G_n(v)|=\frac{\Delta_n(b_n)}{\sqrt{n}b_n}
 $$
where $\Delta_n$ is the modulus of continuity of $G_n$ :
 $
 \Delta_n (b) =\sup_{|u-v|\leq b} |G_n(u)-G_n(v)|.
 $
But we have seen that:
$
||A_nF-f||_\infty = O(b_n).
$
Now, $(b_n)$ satisfies $b_n\to 0$ and $\log n=o(nb_n)$. Set :
$
\kappa_n=(\log n)^{1/2}\log\log n.
$
Then, by Theorem \ref{lemme1}:
$$
\Delta_n(b_n)=O_{a.s.}(\sqrt{b_n\log n})+o_{a.s}(b_n\kappa_n).
$$
Accordingly,
$$
||A_{n}F_n-f||_\infty \leq O_{a.s.}(\frac{\sqrt{\log n}}{\sqrt{n b_n}})+o_{a.s}(\frac{\kappa_n}{\sqrt{n}})+O(b_n)\longrightarrow_{n\to\infty} 0
$$
But we have the following bound:
\begin{eqnarray*}
||B_{n}F_n-f||_\infty & \leq &  2||A_{n}F_n-f||_\infty+||\tau_{b_n/2}f-f||_\infty+||\tau_{-b_n/2}f-f||_\infty
\end{eqnarray*}
The two last terms of the bound are of order $O(b_n)$, hence the result.The optimal rate follows by a direct computation. $\square$

\section{Examples}\label{exemples}

\subsection{ARMA($p$,$q$) models}

These well-known models can be written in the form :
$$
X_n=a_0+\sum_{j=1}^p a_j X_{n-j}+\sum_{j=0}^q b_j\varepsilon_{t-j}
$$
where $a_0,\dots,a_p,b_0,\dots,b_q$ are scalars such that $\sum_{j=0}^q b_j\neq 0$ and $b_0=1$ and where $\{\varepsilon_n,n\in\rel\}$ is a gaussian white noise with distribution $\mathcal{N}(0,\sigma_\varepsilon^2)$.\\
As it is standard in the literature, we assume that the complex polynomials
$A(z)=z^p-\sum_{j=1}^p a_j z^{p-j}$ and $B(z)=\sum_{j=0}^q b_j z^{l-j},\quad (b_0=1)$ have all their roots inside the unit circle.
Then, $\{X_n,n\in\rel\}$ is strictly stationary and its marginals are all gaussian, hence condition $\mathcal{A}_3$ is satisfied. Moreover, one has the representation :
$
X_n=\mu_X+\sum_{j=0}^\infty \beta_j\varepsilon_{n-j}
$
with $\sum_{j=0}^\infty \beta_j^2<\infty$, $\beta_0=1$ and $\mu_X=\esper[X_n]$. Hence, Conditions $\mathcal{A}_1,\mathcal{A}_2$ are satisfied.

\subsection{Linear models}

Consider time series defined by
$
X_n=\mu_X+\sum_{k=0}^\infty a_k\varepsilon_{n-k}
$
where $\mu_X$ is a scalar and where $\{\varepsilon_n,n\in\rel\}$ is a white noise with finite variance. Then, assumption $\mathcal{A}_2$ writes
$
\sum_{k=0}^\infty |a_k|<\infty
$
which is the standard weak dependence condition for linear time series. Moreover, assumption $\mathcal{A}_1$ is satisfied if
$
\sum_{k=0}^\infty a_k^2<\infty.
$

\subsection{Nonlinear autoregressive models}

Consider a strictly stationary Markov process  $\{X_n,n\in\rel\}$ satisfying :
\begin{equation}\label{nlar}
X_{n}=r(X_{n-1})+\varepsilon_n
\end{equation}
where $\{\varepsilon_n,n\in\rel\}$ is a sequence of innovations with absolutely continuous distribution function and finite variance. Let us assume that
$
\rho=\sup_x|r^\prime(x)| <1.
$
Then,
$
Var(X_k-X_k^*) =O(\rho^{2k})
$
and assumption $\mathcal{A}_2$ is satisfied \cite{Diaco99,Wu06}. More generally, assume that 
$
\rho_1(r) < 1.
$
Then, equation (\ref{nlar}) defines a unique stationary distribution \cite{Diaco99} and
$
Var(X_k-X_k^*)=O(\rho_1(r)^{2k}).
$
Assumption $\mathcal{A}_2$ is thus satisfied.\\
Assumption $\mathcal{A}_1$ is satisfied if $r$ is bounded, for instance.\\
Finally, note that tandard examples of time series admitting the representation (\ref{nlar}) are the Threshold Auto-Regressive models of the form :
$$
X_n=a\max(X_{n-1},0)+b\min(X_{n-1},0)+\varepsilon_n
$$ 
where $a$, $b$ are real constants \cite{Wu06}.\\

\noindent {\bf Aknowledgement}. The author would like to thank Dr. B\'{e}atrice Vedel for fruitful discussions which helped simplify the assumptions made in this paper.

\bibliographystyle{elsarticle-num}
\bibliography{mabibfp}

\end{document}